# Euler characteristic of a mixed Koszul complex and the index of vector fields

Achim Hennings[1]

**Abstract**: A three-term complex $R \to R^n \to R^m$ over a local ring determines a mixed Koszul complex, whose Euler characteristic can be expressed by mixed multiplicities. As an application, we offer a simple formula for the index of a holomorphic vector field at an isolated complete intersection singularity.

**0 Introduction:** Let $R$,$\mathfrak{m}$ be a Noetherian local ring of dimension $d$ with its maximal ideal. Let $m$ be a natural number $\geq 1$. Suppose we are given $n = d + m$ elements $a_1, \dots, a_n \in \mathfrak{m}$ and $A_1, \dots, A_n \in \mathfrak{m}R^m$ subject to the relation

(*) $$a_1 A_1 + \dots + a_n A_n = 0.$$

We form the matrix of the columns $A_j$,

$$A = (A_1, \dots, A_n) = (A_{ij}) \in M_{m \times n}(R),$$

and define the submodule $M = (A) \subseteq R^m$ generated by them. Similarly, let $a = (a_1, \dots, a_n)$ and $\mathfrak{a} = (a) \subseteq R$ the ideal generated.

By the relation (*), for any integer $\nu \geq 0$, we can form a double complex $K_\nu$:

$$
K_\nu^{\bullet\bullet} \qquad
\begin{array}{ccccccc|c}
 & & & & & \Lambda^0 \otimes S_{\nu+n} & & n \\
 & & & & & \uparrow & & \\
 & & & & \iddots & & & \\
 & & & & & \vdots & & \\
 & & & \iddots & & & & \\
 & & & & & \uparrow & & \\
 & & \Lambda^0 \otimes S_{\nu+1} \longrightarrow & & \cdots & \longrightarrow \Lambda^{n-1} \otimes S_{\nu+1} & & 1 \\
 & & \uparrow\ d'' & & & \uparrow & & \\
\Lambda^0 \otimes S_\nu & \longrightarrow & \Lambda^1 \otimes S_\nu & \overset{d'}{\longrightarrow} & \cdots & \longrightarrow \Lambda^n \otimes S_\nu & & 0 \\
\hline
0 & & 1 & & & n & &
\end{array}
$$

where we use shorthand notation $\Lambda^\bullet = \Lambda^\bullet R^n$, $S_\bullet = S_\bullet R^m = R[x_1, \dots, x_m]$ and the maps are

$$d' = (\sum_{j=1}^n a_j e_j) \wedge \cdot \otimes id,\ d'' = \sum_{i=1}^m (\sum_{j=1}^n A_{ij} e_j^* \llcorner \cdot) \otimes x_i$$

($\llcorner$ denotes the contraction). The columns are part of the generalized Koszul complexes in the sense of [No, app. C].

We assume $\mathfrak{a}$ and $M$ to be of finite colength. Then the homology module $H(K_\nu)$ of $K_\nu$ as a total complex is of finite length. For the Euler characteristic (alternating sum of the lengths of $H(K_\nu)$ in various degrees) we have:

[1] Universität Siegen, Fakultät IV, Hölderlinstraße 3, D-57068 Siegen

**(0.1) Theorem:** $\chi(K_\nu) = \sum_{i=0}^{d}(-1)^i\, e(\mathfrak{a}, i; M)$.

Here, $e(\mathfrak{a}, i; M)$ $(0 \le i \le d)$ is the $i$th *mixed multiplicity* of the Ideal $\mathfrak{a}$ and the submodule $M \subseteq R^m$. We also define the mixed multiplicity of a *mixed system of parameters*, i.e. of elements $b_1, \dots, b_k \in \mathfrak{m}$ and $B_1, \dots, B_l \in \mathfrak{m}R^m$ with $l = d - k + m - 1$ such that the length of the $R$ −module $R^m/(B) \otimes R/(b)$ is finite: $e(b; B)$. For Cohen-Macaulay rings $e(b; B)$ is simply the length of $R^m/(B) \otimes R/(b)$.

With the assumption that the submodule $(a_1A_1, \dots, a_nA_n) \subseteq R^m$ have finite colength, we get a formula, which is easier to compute.

**(0.2) Theorem:** $\chi(K_\nu) = \sum_{i=0}^{d}(-1)^i\, e(a_1, \dots, a_i; A_{i+2}, \dots, A_n)$.

**(0.3) Remark:** In this article we make essential use of the fact, that a complete generalized Koszul complex (in the sense of [No, app. C]) associated to an $m \times n$ −matrix with $n > d + m - 1$ has Euler characteristic zero ([BV], [Ki]). In particular this applies to the ordinary Koszul complex for $> d$ elements. Using this fact, we may modify the complex by adding or canceling of rows or columns with $\chi = 0$ without changing the Euler characteristic of the total complex. In this way, we can easily see that $\chi(K_\nu)$ is independent of $\nu$.

**(0.4) Remark:** The relation (*) means, that there is a complex $R \to R^n \to R^m$. For a complex $R^r \to R^n \to R^m$ the analog would be generalized Koszul complexes in both directions and mixed multiplicities of two modules. Cf. [W] for the most general construction.

**(0.5) Remark:** By deleting the last $m$ columns of $K_0^{\bullet\bullet}$ (having $\chi = 0$) and completing the Koszul complexes in the remaining columns we obtain the new complex $L^{\bullet\bullet}$ below. By remark 0.3, $\chi(L^{\bullet\bullet}) = 0$. Let $L_+^{\bullet\bullet}$ resp. $L_-^{\bullet\bullet}$ denote the part of $L^{\bullet\bullet}$ above resp. below the dotted line. Then $\chi(L^{\bullet\bullet}) = \chi(L_+^{\bullet\bullet}) + \chi(L_-^{\bullet\bullet}) = 0$ and $\chi(L_+^{\bullet\bullet}) = \chi(K_0^{\bullet\bullet})$.

$$
\begin{array}{cccccccccc}
 & 0 & & 1 & & & & n-m & & \\
 & & & & & & & \Lambda^0\otimes S_{n-m} & & n \\
 & & & & & & & \uparrow & & \\
 & & & & & \cdot^{\cdot^{\cdot}} & & \vdots & & \\
L^{\bullet\bullet} & & & & \cdot^{\cdot^{\cdot}} & & & \uparrow & & \\
 & & & \Lambda^0\otimes S_1 & \longrightarrow & \cdots & \longrightarrow & \Lambda^{n-m-1}\otimes S_1 & & 1 \\
 & & & \uparrow\ d'' & & & & \uparrow & & \\
 & \Lambda^0\otimes S_0 & \longrightarrow & \Lambda^1\otimes S_0 & \overset{d'}{\longrightarrow}\cdots & \Lambda^{n-m-1}\otimes S_0 & \longrightarrow & \Lambda^{n-m}\otimes S_0 & & 0 \\
 & \uparrow & & \uparrow & & \uparrow & & \uparrow & & \\
\hline
 & \Lambda^m\otimes S_0 & \longrightarrow & \Lambda^{m+1}\otimes S_0 & \longrightarrow\ \cdots & \Lambda^{n-1}\otimes S_0 & \longrightarrow & \Lambda^n\otimes S_0 & & -1 \\
 & \uparrow & & \uparrow & & \uparrow & & & & \\
 & \Lambda^{m+1}\otimes S_1 & \longrightarrow & \Lambda^{m+2}\otimes S_1 & \longrightarrow\ \cdots & \Lambda^n\otimes S_1 & & & & \\
 & \vdots & & & \cdot^{\cdot^{\cdot}} & & & & & \\
 & \uparrow & & & & & & & & \\
 & \Lambda^n\otimes S_{n-m} & & \cdot^{\cdot^{\cdot}} & & & & & & -n-m-1
\end{array}
$$

**(0.6) Application: Index of a vector field at a complete intersection singularity**

Let $(X,0) \subseteq (\mathbb{C}^n,0)$ be an analytic isolated complete intersection singularity of dimension $d$, defined by a regular sequence $f=(f_1,\dots,f_m)$. Let $V=\sum_{i=1}^n V_i \frac{\partial}{\partial x_i}$ be a holomorphic vector field on $X$ with an isolated zero. Then we have the relation $\sum_{i=1}^n V_i \frac{\partial f}{\partial x_i}=0$ in $\mathcal{O}_X^m$ and we are in the situation as above.

The vector field $V$ also defines a contraction $i_V\colon \Omega_X^p \to \Omega_X^{p-1}$ and, since $(i_V)^2=0$, a complex

$$(\Omega_{X,0}^{\bullet}, i_V) = (\Omega_{X,0}^d \to \cdots \to \Omega_{X,0}^0).$$

We consider the Euler characteristic of this complex,

$$\chi(\Omega_{X,0}^{\bullet}, i_V),$$

viewed as a homological complex in degrees $d,\dots,0$

A natural definition of the *index* of $V$ at the isolated zero, generalizing the smooth case to isolated complete intersections, is given by

$$ind_X(V,0) = \chi(\Omega_{X,0}^{\bullet}, i_V) - (-1)^d \mu(X,0),$$

where $\mu$ is the Milnor number. This definition satisfies the Hopf index formula and the usual continuity under deformation. (Cf. [BEG] and the references therein.)

We form the complex $L^{\bullet\bullet}$ with $R=\mathcal{O}_{X,0}$. By the CM property, the columns of $L_{-}^{\bullet\bullet}$ are exact, except at the end. Then it is easy to see, using $\Lambda^p R^n \cong \Lambda^{n-p} R^n$, that $L_{-}^{\bullet\bullet}$ is a resolution of $(\Omega_{X,0}^{-\bullet+d}, i_V)$. We obtain from remark 0.5

$$\chi(\Omega_{X,0}^{\bullet}, i_V) = -(-1)^d \chi(L_{-}^{\bullet\bullet}) = (-1)^d \chi(K_0^{\bullet\bullet})$$

By theorem 0.2 we thus obtain a simple length formula for the index $ind(V,0)$:

$$ind_X(V,0) = (-1)^d \chi(K_0^{\bullet\bullet}) - (-1)^d \mu(X,0),$$

$$\chi(K_0^{\bullet\bullet}) = \textstyle\sum_{i=0}^d (-1)^i e(V_1,\dots,V_i; \frac{\partial f}{\partial x_{i+2}},\dots,\frac{\partial f}{\partial x_n}).$$

Here we have to choose coordinates such that $(V_1 \frac{\partial f}{\partial x_1},\dots,V_n \frac{\partial f}{\partial x_n})$ has finite colength (cf. Lemma 3.1.3).

## 1 Mixed multiplicities of a module and an ideal

**1.1** We give a short exposition of the known theory of mixed multiplicities in the easy case of a submodule and an ideal of finite colength.

Let $R$,$\mathfrak{m}$ a Noetherian local ring of dimension $d$. Let be given an ideal $\mathfrak{a} \subseteq \mathfrak{m}$ with generators $a_1,\dots,a_k$ and a submodule $M \subseteq \mathfrak{m}R^m$ with its generating matrix $A=(A_1,\dots,A_n)$, both of finite colength. Then $n \geq m+d-1$ if $d>0$. We consider the symmetric algebras

$$S = SR^m = R[x_1,\dots,x_m],\ S(A) = SM$$

and the image of $S(A)$ in $S$,

$$U = im\left(S(A)\right) = R[y_1, \dots, y_n] \subseteq S,\ y_j := \sum_{i=1}^{m} A_{ij}\, x_i.$$

**(1.1.1) Lemma:** For $\mu, \nu$ >>0, the length $L(S_\nu/\mathfrak{a}^\mu U_\nu)$ is a rational polynomial $Q(\mu,\nu)$ of degree $d+m-1$ with leading term representable as

$$Q^* = \sum_{i=0}^{d} \frac{e_i}{i!(m+d-1-i)!}\mu^i \nu^{m+d-1-i},\ e_i \in \mathbb{N}^+.$$

The coefficient $e_i = e(\mathfrak{a}, i; M)$ is called the $i$th *mixed multiplicity* of $\mathfrak{a}, M$. In particular, $e_d = e(\mathfrak{a})$ is the multiplicity of the ideal and $e_0 = e(M)$ is the multiplicity of the submodule in the sense of Buchsbaum and Rim [BR] (denoted $e(R^m/M)$ there). More generally, $e(\mathfrak{a}, i; M; E)$ can be defined for finitely generated modules $E$ if $R/\mathfrak{a}\otimes E$ and $R^m/M\otimes E$ are of finite length. In this case $d = \dim E$ replaces $\dim R$.

Proof: 1) By [BR] or [KT] the length $L(S_\nu/U_\nu) = P(\nu)$ for $\nu \gg 0$ is a rational polynomial of degree $d+m-1$, and $e(M)$ is the $(d+m-1)$!-fold of the leading coefficient.

2) The ring

$$\bar{U} := \bigoplus_{\mu,\nu\geq 0} \mathfrak{a}^\mu U_\nu/\mathfrak{a}^{\mu+1}U_\nu$$

is a homogeneous quotient of the bi-graded ring $R/\mathfrak{a}[X_1, \dots, X_k; Y_1, \dots, Y_n]$ (mapping $X_i \mapsto \bar{a}_i$, $Y_j \mapsto \bar{A}_j$). Therefore, $L\left(\bar{U}_{\mu\nu}\right) = q(\mu,\nu)$ is a rational polynomial for $\mu, \nu$ >>0. Similarly $L(U_\nu/\mathfrak{a}^{\mu_0}U_\nu)$ for $\nu \gg 0$ and $\mu_0$ fixed is a rational polynomial $q_{\mu_0}(\nu)$ whose degree is $< d + m - 1$ by Lemma (4.1.3) of the appendix.

3) Thereby

$$L(S_\nu/\mathfrak{a}^\mu U_\nu) = P(\nu) + q_{\mu_0}(\nu) + \sum_{i=\mu_0}^{\mu-1} q(i,\nu)$$

for $\mu_0$ sufficiently large and $\mu, \nu$ >>0 is a rational polynomial $Q(\mu,\nu)$ as well. The degree is $d+m-1$, since $Q(\mu_0\nu, \nu)$ for any $\mu_0 \gg 0$ and fixed is the Buchsbaum-Rim polynomial of the submodule $\mathfrak{a}^{\mu_0}M$ and of degree $d+m-1$. The degree in $\mu$ is the degree of $L\left(S_{\nu_0}/\mathfrak{a}^\mu U_{\nu_0}\right)$ for $\nu_0$ fixed and large. Since $L(S_{\nu_0}/U_{\nu_0}) < \infty$, this polynomial has even the same leading term as $L\left(S_{\nu_0}/\mathfrak{a}^\mu S_{\nu_0}\right)$, namely

$$\frac{e(\mathfrak{a})}{d!}\binom{\nu_0 + m - 1}{m-1}\mu^d.$$

This shows $e_d = e(\mathfrak{a})$. The representation of the leading term of $Q$ as

$$Q^* = \sum_{i=0}^{d} \frac{e_i}{i!(m+d-1-i)!}\mu^i \nu^{m+d-1-i},\ e_i \in \mathbb{Z},$$

now follows from the fact, that $Q$ is numerical, i.e. has integer values for $\mu, \nu$ >>0.

By the representation above, $P(\nu)$ and $Q(\mu_0, \nu)$ have the same leading term, showing that $e_0 = e(M)$.

4) There remains to show that the integers $e_i = e(\mathfrak{a}, i; M)$ are positive. We reduce this to the familiar cases $e_0 = e(M) > 0$ and $e_d = e(\mathfrak{a}) > 0$. We need some definitions and facts about superficial elements, which are collected in appendix 4.2. We may assume that $R$ is an integral domain (Lemma 1.1.3) and $|R/\mathfrak{m}| = \infty$. Then there is a superficial non-zero-divisor $a \in \mathfrak{a}$ for $\mathfrak{a}, M$, and we have for $\bar{R} = R/aR$ and the images $\bar{\mathfrak{a}}$, $\bar{M}$

$$e(\mathfrak{a}, i; M) = e(\bar{\mathfrak{a}}, i-1; \bar{M}),\ 1 \le i \le d,\ d \ge 2,$$

and $\dim \bar{R} = d - 1$.

5) The general case of an additional module $E$ can be treated the same way. For the positivity of the mixed multiplicities we again use Lemma 1.1.3.

**(1.1.2) Lemma:** If we define $e(\mathfrak{a}, i; M; E, d) = e(\mathfrak{a}, i; M; E)$ if $\dim E = d$ and $e(\mathfrak{a}, i; M; E, d) = 0$ if $\dim E < d$, we obtain an additive function on modules of dimension $\le d$.

Proof: For a given exact sequence $0 \to E' \to E \to E'' \to 0$ with $\dim E \le d$, the sequence

(*) $\quad 0 \to S_\nu \otimes E'/(\mathfrak{a}^\mu U_\nu E) \cap (S_\nu \otimes E') \to S_\nu \otimes E/\mathfrak{a}^\mu U_\nu E \to S_\nu \otimes E''/\mathfrak{a}^\mu U_\nu E'' \to 0$

is exact as well (where $U_\nu E, \ldots$ denotes the image of $U_\nu \otimes E, \ldots$ in $S_\nu \otimes E, \ldots$).

From the fact, that

$$\bigoplus_{\mu,\nu \ge 0} (\mathfrak{a}^\mu U_\nu E) \cap (S_\nu \otimes E') \subseteq \bigoplus_{\mu,\nu \ge 0} \mathfrak{a}^\mu U_\nu E$$

is a (finitely generated) $(\bigoplus_{\mu,\nu \ge 0} \mathfrak{a}^\mu U_\nu)$-submodule we get $\mu_0, \nu_0$ such that for $\mu \ge \mu_0, \nu \ge \nu_0$:

$$(\mathfrak{a}^\mu U_\nu E) \cap (S_\nu \otimes E') \subseteq \mathfrak{a}^{\mu-\mu_0} U_{\nu-\nu_0}(\mathfrak{a}^{\mu_0} U_{\nu_0} E \cap S_{\nu_0} \otimes E') \subseteq \mathfrak{a}^{\mu-\mu_0} U_{\nu-\nu_0}(S_{\nu_0} \otimes E')$$

Let $\mu_1$ be large enough, such that $\mathfrak{a}^{\mu_1} S_{\nu_0} \subseteq U_{\nu_0}$. Then for $\mu \ge \mu_0 + \mu_1$:

$$\mathfrak{a}^\mu U_\nu E' \subseteq (\mathfrak{a}^\mu U_\nu E)) \cap (S_\nu \otimes E') \subseteq \mathfrak{a}^{\mu-\mu_0-\mu_1} U_\nu E'.$$

Hence the $\mu^i \nu^{m+d-1-i}$-term of $L(S_\nu \otimes E'/(\mathfrak{a}^\mu U_\nu E) \cap (S_\nu \otimes E'))$ (which is a polynomial by (*)) just provides $e(\mathfrak{a}, i; M; E', d)$.

**(1.1.3) Lemma:** With $\mathfrak{p}_1, \ldots, \mathfrak{p}_r$ the minimal prime ideals of $E$ of maximal dimension $d = \dim E$, we have

$$e(\mathfrak{a}, k; M; E) = \textstyle\sum_{i=1}^r e(\mathfrak{a}, k; M; R/\mathfrak{p}_i) L(E_{\mathfrak{p}_i}).$$

Proof: We consider a chain $E = E_0 \supseteq \cdots \supseteq E_s = 0$ of submodules such that $E_j/E_{j+1} \cong R/\mathfrak{p}(j)$ is a quotient by an associated prime ideal. Then $e(\mathfrak{a}, k; M; E)$ is the sum of all $e(\mathfrak{a}, k; M; R/\mathfrak{p}(j), d)$. The number of occurrences of a minimal prime ideal $\mathfrak{p}_i$ among the $\mathfrak{p}(j)$ is just $L(R_{\mathfrak{p}_i})$.

## 1.2 Multiplicity of parameter systems

Let $R$,$\mathfrak{m}$ be a Noetherian local ring, $a \in \mathfrak{m}R^k$ a row-vector and $A = (A_1, \ldots, A_n) \in \mathfrak{m}M_{m \times n}(R)$ a matrix. We consider the associated Koszul complex,

$$K: 0 \to \Lambda^0 R^k \to \cdots \to \Lambda^k R^k \to 0$$

(with index range $-k, \ldots, 0$), and the generalized Koszul complex ([No, appendix C]),

$$L_\nu: 0 \to S_\nu R^m \otimes \Lambda^0 R^n \to \cdots \to S_{\nu+n} R^m \otimes \Lambda^n R^n \to 0,\ \nu \ge 0$$

(with index range $-n, \ldots, 0$).

In the following let $E$ be a finitely generated $R$-module of dimension $d$, such that the length of $R/(a)\otimes R^m/(A)\otimes E$ is finite[2]. If $d-k>0$, the inequality $n\geq d-k+m-1$ must hold.

In the equality case we call $a, A$ a *mixed system of parameters*. This is the essential case of the following

**(1.2.1) Definition:** $e(a;A;E) \coloneqq \chi(K\otimes L_\nu\otimes E)$.

**(1.2.2) Lemma:** $\chi(K\otimes L_\nu\otimes E)$ is independent of $\nu\geq 0$.

Proof: By hypothesis, the homology of the $E_1$-term $L_\nu\otimes H^\bullet(K\otimes E)$ is of finite length. Therefore $\chi(K\otimes L_\nu\otimes E)=\chi(L_\nu\otimes H^\bullet(K\otimes E))=\sum_i(-1)^i\chi(L_\nu\otimes H^i(K\otimes E))$, and we use the fact from [Ki], that $\chi(L_\nu\otimes G)$ for a finitely generated module $G$ is independent of $\nu$.

We remark here, that $e(a,A;E)$ depends on $A$ only modulo $(a)$. Namely, the module $H^\bullet(K\otimes E)$ in the foregoing argument is annihilated by $(a)$.

**(1.2.3) Lemma:** $e(a;A;E)=0$ in the cases $n>d-k+m-1$ or $k-d>0$.

Proof: In the case $k-d>0$ the assertion follows from $\chi(H^\bullet(L_\nu\otimes E)\otimes K)=0$.

Let us assume $n>d-k+m-1$. From the additivity of the Euler characteristic we obtain as in Lemma 1.1.3 and Lemma 1.2.5 below:

$$\begin{aligned} e(a;A;E) &= \chi(L_\nu\otimes K\otimes E) = \chi\big(L_\nu\otimes H^\bullet(K\otimes E)\big) \\ &= \sum_{\mathfrak{p}\in Min(H(K\otimes E))} \chi(L_\nu\otimes R/\mathfrak{p})\chi\big(H^\bullet(K\otimes E)_\mathfrak{p}\big) \\ &= \sum_{\mathfrak{p}\in Min(H(K\otimes E))} \chi(L_\nu\otimes R/\mathfrak{p})\chi\big(K\otimes E_\mathfrak{p}\big) \end{aligned}$$

In the last sum $\chi\big(K\otimes E_\mathfrak{p}\big)=0$ if $k>\dim E_\mathfrak{p}$, and $\chi(L_\nu\otimes R/\mathfrak{p})=0$ if $n-m+1>\dim R/\mathfrak{p}$ (also when $\dim R/\mathfrak{p}=0$) ([BV], [Ki]), and therefore we need to sum only over minimal primes. Because of

$$\begin{aligned} \dim E_\mathfrak{p}+\dim R/\mathfrak{p} = ht(\mathfrak{p}/Ann(E))+\dim R/\mathfrak{p} &\leq \dim E \\ &< n-m+1+k \end{aligned}$$

at least one of these conditions holds.

**(1.2.4) Lemma:** For $k=d=\dim E$ and $n=(d-k)+m-1=m-1$ we obtain the ordinary multiplicity of the ideal: $e(a;A;E)=e(a;E)$.

**Proof:** The assumption $n=m-1$ implies $\dim E/(a)E=0$ and

$$e(a;A;E)=\chi(H^\bullet(K\otimes E)\otimes L_\nu).$$

The complex $L_\nu$ with $\nu=0$ has alternating rank sum

$$(-1)^n\sum_{i=0}^n(-1)^i\binom{n}{i}\binom{n+i}{n}.$$

This is the $n$th coefficient of the power series

[2] We remark $supp(M_1\otimes M_2)=supp(M_1)\cap supp(M_2)$ for finitely generated modules.

$$\frac{1}{n!}\left(\frac{1}{1-x}\right)^{(n)}(1-x)^n = \frac{1}{1-x},$$

hence equal to 1. We obtain

$$\chi(H^\bullet(K \otimes E) \otimes L_\nu) = e(a; E).$$

In the other extremal case $k = 0$ and $n = m + d - 1$, where $d = \dim E$, $e(a; A; E) = e((A); E)$ is the Buchsbaum-Rim multiplicity of the submodule $(A) \subseteq R^m$ with respect to the module $E$ (cf. [Ki, th. 4]).

**(1.2.5) Lemma:** The multiplicity $e(a; A; E)$ is additive in $E$.

Proof: Let $0 \to E' \to E \to E'' \to 0$ be an exact sequence. Then $e(a; A; E)$ is defined if and only if $e(a; A; R/\mathfrak{p})$ is defined for all $\mathfrak{p} \in supp(E) = supp(E') \cup supp(E'')$, hence if and only if $e(a; A; E')$ and $e(a; A; E'')$ are both defined. From the exact sequence we get an exact sequence of complexes by applying $K \otimes L_\nu \otimes$.

**(1.2.6) Lemma:** With $\mathfrak{p}_1, \dots, \mathfrak{p}_r$ the minimal primes of $E$ of maximal dimension $\dim E$ we have

$$e(a; A; E) = \sum_{i=1}^{r} e(a; A; R/\mathfrak{p}_i) L(E_{\mathfrak{p}_i}).$$

Proof: We apply Lemma 1.2.5 as in Lemma 1.1.3. By Lemma 1.2.3, all summands not appearing are zero.

**(1.2.7) Lemma:** Let $\mathfrak{a} \subseteq \mathfrak{m}$ and $M \subseteq \mathfrak{m}R^m$ of finite codimension. Let $a$ and $A$ as above with $k \leq d = \dim E$, $n = d - k + m - 1$. Assume $(a) \subseteq \mathfrak{a}$, $(A) \subseteq M$. Then

$$e(a; A; E) \geq e(\mathfrak{a}, k; M; E, d).$$

Proof: By Lemma 1.1.3 and Lemma 1.2.6 we may assume that $E = R$ is an integral domain. The case $k = 0$ reduces to ordinary Buchsbaum-Rim multiplicities: $e(a; A) = e\big((A)\big) \geq e(M) = e(\mathfrak{a}, 0; M)$.

In the case $k > 0$ we have $a_1 \neq 0$ and we pass to $\bar{R} = R/a_1R$. By Lemma 1.2.8 below

$$e(a; A; R) = e(a_2, \dots, a_k; A; \bar{R}) - e\big(a_2, \dots, a_k; A; (0{:}\, a_1R)\big) = e(a_2, \dots, a_k; A; \bar{R}).$$

In general (cf. appendix, remark to Lemma 4.2.2) $e(\mathfrak{a}, k; M; R) \leq e(\mathfrak{a}, k-1; M; \bar{R})$. Using induction on $k$ we may assume $e(\mathfrak{a}, k-1; M; \bar{R}) \leq e(a_2, \dots, a_k; A; \bar{R})$. The induction step follows.

**(1.2.8) Lemma:** We have

$$e(a; A; E) = e(a_2, \dots, a_k; A; E/a_1E) - e(a_2, \dots, a_k; A; (0{:}\, a_1)_E).$$

Proof: Let $K'$ denote the Koszul complex of $a_2, \dots, a_k$. Then $K$ is the total complex of $K' \xrightarrow{a_1} K'$ and $K \otimes L_\nu \otimes E$ is the total complex of $K' \otimes L_\nu \otimes E \xrightarrow{a_1} K' \otimes L_\nu \otimes E$. The $E_1$-terms of this double complex are $K' \otimes L_\nu \otimes (0{:}\, a)_E$ and $K' \otimes L_\nu \otimes E/a_1E$. Therefore

$$\chi(K \otimes L_\nu \otimes E) = \chi(K' \otimes L_\nu \otimes E/a_1E) - \chi(K' \otimes L_\nu \otimes (0{:}\, a)_E).$$

**(1.2.9) Lemma:** (additivity)

If one of the two sides is defined, the following formula holds:

$$e(a_1, \dots, a_{k-1}; a_k A_k, A_{k+1}, \dots, A_n; E)$$
$$= e(a_1, \dots, a_k; A_{k+1}, \dots, A_n; E) + e(a_1, \dots, a_{k-1}; A_k, A_{k+1}, \dots, A_n; E)$$

Proof: By Lemma (1.2.8) we can reduce to the case $k = 1$. Let $L'_\nu$ be the Koszul complex of $A_{k+1}, \dots, A_n$ and $L'_\nu \xrightarrow{A_k} L'_{\nu+1}$ the multiplication map by $\sum_{i=1}^m A_{ik}\, x_i$. This represents $L_\nu$ as a double complex. We consider the composition

$$L'_\nu \xrightarrow{a_k} L'_\nu \xrightarrow{A_k} L'_{\nu+1},$$

and we have to show:

$$\chi\left(L'_\nu \xrightarrow{a_k} L'_\nu\right) + \chi\left(L'_\nu \xrightarrow{A_k} L'_{\nu+1}\right) = \chi\left(L'_\nu \xrightarrow{a_k A_k} L'_{\nu+1}\right).$$

The Euler characteristic of these complexes can be computed via the $E_2$-terms, the $E_1$-terms being the homology of $L'_\nu$ and $L'_{\nu+1}$. We have the diagram

$$(*) \qquad H(L'_\nu) \xrightarrow{H(a_k)} H(\, L'_\nu) \xrightarrow{H(A_k)} H(L'_{\nu+1})$$

and we need to compare the lengths of the kernel and cokernel of $H(a_k)$, $H(A_k)$, $H(a_k A_k)$.

For a composition of module homomorphisms $M_1 \xrightarrow{f} M_2 \xrightarrow{g} M_3$ there is an exact kernel-cokernel-sequence

$$0 \to \ker f \to \ker gf \to \ker g \to \operatorname{coker} f \to \operatorname{coker} gf \to \operatorname{coker} g \to 0.$$

If we apply this to (*) (in each degree separately) we get as we need

$$L(\operatorname{coker} H(a_k)) - L(\ker H(a_k)) + L(\operatorname{coker} H(A_k)) - L(\ker H(A_k)) =$$
$$L(\operatorname{coker} H(a_k A_k)) - L(\ker H(a_k A_k)).$$

**(1.2.10) Lemma:** When $E$ is a Cohen-Macaulay module of dimension $d > k$ and $n = d - k + m - 1$, the mixed multiplicity is a length,

$$e(a; A; E) = L(R/(a) \otimes R^m/(A) \otimes E).$$

Proof: We have $H(K \otimes E) = H^0(K \otimes E) = E/(a)$ since $a$ is a regular sequence. Also the generalized Koszul complex

$$0 \to S_{n-m-1} R^m \otimes \Lambda^n R^n \otimes M \to \cdots \to S_0 R^m \otimes \Lambda^{m+1} R^n \otimes M \to S_0 R^m \otimes \Lambda^1 R^n \otimes M$$
$$\to S_1 R^m \otimes \Lambda^0 R^n \otimes M \to 0$$

with $M = E/(a)E$ has its homology $E/(a)E \otimes R^m/(A)$ only in dimension zero (cf. [No], [Ki]). By [Ki], the Euler characteristic is the same as of the complex $L_\nu \otimes M$ (of modified degree). We see: $\chi(L_\nu \otimes K \otimes E) = \chi(L_\nu \otimes E/(a)) = L(E/(a)E \otimes R^m/(A))$.

**(1.2.11) Lemma:** For $n = d - k + m - 1$ and $d \geq k$,

$$e(a; A; E) = \textstyle\sum_{\mathfrak{p}} e(a; E_{\mathfrak{p}}) e(A; R/\mathfrak{p}),$$

where $\mathfrak{p}$ runs over the minimal primes containing $(a) + Ann(E)$.

(A more general associativity formula can be formulated for part of the sequence $a$.)

Proof: Similar to the proof of Lemma 1.2.3 and Lemma 1.2.6,

$$\chi(L_\nu \otimes K \otimes E) = \chi(L_\nu \otimes H^\bullet(K \otimes E)) = \sum_{\mathfrak{p}} \chi(L_\nu \otimes R/\mathfrak{p}) \chi\big(H^\bullet(K \otimes E)_{\mathfrak{p}}\big).$$

## 2 Proof of Theorem (0.1)

For the given matrix $A = (A_1, \dots, A_n) \in M_{m \times n}$ we already have put $y_j = \sum_{i=1}^m A_{ij} x_i$ and $U = R[y_1, \dots, y_n] \subseteq S = R[x_1, \dots, x_m]$. We define a filtration on $K_\nu$ by subcomplexes $F^\mu K_\nu$,

$$F^\mu K_\nu^{pq} = \mathfrak{a}^{\mu+p} \Lambda^{p-q} U_{\nu+q} \subseteq \Lambda^{p-q} \otimes S_{\nu+q} = K_\nu^{pq}.$$

We determine $\chi(K_\nu / F^\mu K_\nu)$ by an adaptation of the usual proof for the ordinary Koszul complex (cf. [Se]).

**(2.1) Lemma:** For $\mu, \nu \gg 0$ the complex $F^\mu K_\nu$ is acyclic.

Proof: We define the bi-graded ring

$$\widehat{U} = \oplus_{\mu,\nu \geq 0} \mathfrak{a}^\mu U_\nu s^\mu \subseteq R[x][s].$$

With $G = \widehat{U}^n$ we form the double complex:

$$\begin{array}{lccccccl}
 & & & & & & \Lambda^0 G & \\
 & & & & & \iddots & \vdots & \\
\widehat{K}: & & & & & & & \\
 & & & \iddots & & & \uparrow & \\
 & & \Lambda^0 G & \longrightarrow & \cdots & \longrightarrow & \Lambda^{n-1} G & \\
 & & d'' \uparrow & & & & \uparrow & d' = m_a \\
 \Lambda^0 G & \longrightarrow & \Lambda^1 G & \overset{d'}{\longrightarrow} & \cdots & \longrightarrow & \Lambda^n G & d'' = i_A
\end{array}$$

The columns are ordinary Koszul complexes defined by contraction $i_A$ with the linear form $A = \sum_{j=1}^n y_j e_j^*$ and the rows are defined by exterior multiplication $m_a$ with $a = \sum_{j=1}^n a_j s e_j \in G_{10}$. We define the subcomplex of $\widehat{K}$,

$$FK := \oplus_{\mu,\nu \geq 0} F^\mu K_\nu.$$

This is a finitely generated module over $\widehat{U}$. For any $j, k \in \{1, \dots, n\}$ we define a homotopy operator $h = h' + h''$ on $FK$ with horizontal and vertical components as follows:

$$h'|_{pq} = y_k i_{e_j^*},$$

$$h''|_{pq} = \begin{cases} m_{e_k} i_{e_j^*} m_a, \text{if } p = n \\ 0, \text{else} \end{cases}$$

(where the operators act on $\Lambda^\bullet G$, and the result is taken as an element of the appropriate bi-degree).

Assertion 1: $dh + hd = (a_j s) y_k$

Proof: The left-hand side splits into four summands:

$$d''h' + h'd'' = i_A y_k i_{e_j^*} + y_k i_{e_j^*} i_A = 0$$

$$d'h' + h'd'|_{pq} = m_a y_k i_{e_j^*} + y_k i_{e_j^*} m_a = y_k \left( m_a i_{e_j^*} + i_{e_j^*} m_a \right) = y_k (a_j s) \text{ if } p \neq n,$$

$$d'h' + h'd'|_{pq} = y_k m_a i_{e_j^*}, \text{ if } p = n,$$

$$d'h'' + h''d'|_{pq} = 0, \text{ if } p \neq n-1,$$

$$d'h'' + h''d'|_{pq} = m_{e_k} i_{e_j^*} m_a m_a = 0, \text{ if } p = n-1,$$

$$d''h'' + h''d''|_{pq} = 0, \text{ if } p \neq n,$$

$$\begin{aligned} d''h'' + h''d''|_{pq} &= i_A m_{e_k} i_{e_j^*} m_a + m_{e_k} i_{e_j^*} m_a i_A \\ &= y_k i_{e_j^*} m_a - m_{e_k} i_A i_{e_j^*} m_a + m_{e_k} i_{e_j^*} m_a i_A \\ &= y_k i_{e_j^*} m_a + m_{e_k} i_{e_j^*} (i_A m_a + m_a i_A) = y_k i_{e_j^*} m_a \end{aligned}$$

(since $\sum_j y_j a_j = 0$), if $p = n$.

The sum of all the four terms is in any case $y_k(a_j s)$.

From assertion 1 we see, that the homology module $H(FK)$ is annihilated by $\widehat{U}_{11} = \mathfrak{a}(y_1, \dots, y_n)s$. Therefore, all homogeneous components $(\mu, \nu) \geq (\mu_0 + 1, \nu_0 + 1)$ are zero if $(\mu_0, \nu_0)$ is an upper bound for the degrees of generators of $H(FK)$. This was to be shown.

**(2.2) Lemma:** For $\mu, \nu \gg 0$, $\chi(K_\nu / F^\mu K_\nu) = \sum_{i=0}^{d} (-1)^i e_i$, where $e_i = e(\mathfrak{a}, i; M)$.

Proof: Making use of the polynomial $Q(\mu, \nu) = L(S_\nu / \mathfrak{a}^\mu U_\nu)$ from Lemma 1.1.1 we compute (with $k = p - q$)

$$\chi(K_\nu / F^\mu K_\nu) = \sum_{p=0}^{n} \sum_{q=0}^{p} (-1)^{p+q} \binom{n}{p-q} Q(\mu + p, \nu + q) =$$
$$\sum_{k=0}^{n} (-1)^k \binom{n}{k} \left( \sum_{p=k}^{n} Q(\mu + p, \nu + p - k) \right) = (-1)^n \Delta_k^n|_0 \left( \sum_{p=k}^{n} Q(\mu + p, \nu + p - k) \right) =$$
$$\sum_{i=0}^{d} (-1)^i e_i.$$

The last equality holds by the following

**(2.3) Lemma:** For $a + b = n - 1$,

$$\Delta_k^n|_0 \left( \sum_{p=k}^{n} (\mu + p)^a (\nu + p - k)^b \right) = -\Delta_k^n|_0 \left( \sum_{p=0}^{k-1} (\mu + p)^a (\nu + p - k)^b \right) = -(-1)^b a!\, b!.$$

For $a + b < n - 1$ the result is zero.

Proof: The first equality is true because $\sum_{p=0}^{n} (\mu + p)^a (\nu + p - k)^b$ is of degree $\leq n - 1$ in $k$. The second equality we prove by induction on $b$. For $b = 0$ we have

$$\Delta_k^n|_0 \left( \sum_{p=0}^{k-1} (\mu + p)^a \right) = \Delta_k^{n-1}|_0 (\mu + k)^a = \begin{cases} (n-1)!, a = n - 1 \\ 0, a < n - 1 \end{cases}.$$

For $b > 0$ we have

$$\Delta_k^n|_0\left(\sum_{p=0}^{k-1}(\mu+p)^a(\nu+p-k)^b\right) = \Delta_k^{n-1}|_0\left(\sum_{p=0}^{k}(\mu+p)^a(\nu+p-k-1)^b - \sum_{p=0}^{k-1}(\mu+p)^a(\nu+p-k)^b\right) = \Delta_k^{n-1}|_0\left((\mu+k)^a(\nu+k-k-1)^b + \sum_{p=0}^{k-1}(\mu+p)^a((\nu+p-k-1)^b-(\nu+p-k)^b)\right).$$

Here $(\nu+p-k-1)^b-(\nu+p-k)^b = (-1)b(\nu+p-k)^{b-1}+q(\nu+p-k)$ with $q$ of degree $\leq b-2$ and $(\mu+k)^a(\nu-1)^b$ is of degree $a \leq n-2$ in $k$. We thus can apply the induction hypothesis to conclude.

Proof of Theorem 0.1:

By Remark 0.3, $\chi(K_\nu)$ is independent of $\nu$. For $\mu,\nu \gg 0$, $\chi(K_\nu) = \chi(K_\nu/F^\mu K_\nu) = \sum_{i=0}^{d}(-1)^i e_i$ by Lemma 2.1 and Lemma 2.2.

## 3 Proof of Theorem 0.2

### 3.1 Alternating multiplicity

Let $R,\mathfrak{m}$ be a Noetherian local ring of dimension $d$. For elements $a_1,\dots,a_n \in \mathfrak{m}$ and $A_1,\dots,A_n \in \mathfrak{m}R^m$ with $n = d+m$ may hold the relation

$$a_1A_1 + \dots + a_nA_n = 0,$$

and the submodule $(a_1A_1,\dots,a_nA_n)$ be of finite codimension.

In this case, we say that the *alternating multiplicities* are defined, which are given by the following formula:

For each injective sequence $\boldsymbol{i} = (i_1,\dots,i_{d+1})$ from $\{1,\dots,n\}$ we form the alternating sum of multiplicities

$$S(a,A,\boldsymbol{i}) = \sum_{j=0}^{d}(-1)^j\, e(a_{i_1},\dots,a_{i_j}; A_k, k \notin \{i_1,\dots,i_{j+1}\}).$$

By the following result, we henceforth may simply write $S(a,A,\boldsymbol{i}) = S(a,A)$.

**(3.1.1) Lemma:** $S(a,A,\boldsymbol{i})$ is independent of $\boldsymbol{i}$.

Proof: In this proof we put $S(\boldsymbol{i}) = S(a,A,\boldsymbol{i})$.

Assertion 1: $S(i_1,\dots,i_d,i) = S(i_1,\dots,i_d,j)$

This holds, since all summands are equal with possible exception of

$$(-1)^d e(a_{i_1},\dots,a_{i_d}; A_k, k \notin \{i_1,\dots,i_d,i\}) \text{ and } (-1)^d e(a_{i_1},\dots,a_{i_d}; A_k, k \notin \{i_1,\dots,i_d,j\}).$$

However, they have the same value $e(a_{i_1},\dots,a_{i_d})$ by Lemma 1.2.4.

Assertion 2: $S(i_1,\dots,i_{k-1},i,j,i_{k+2},\dots,i_{d+1}) = S(i_1,\dots,i_{k-1},j,i,i_{k+2},\dots,i_{d+1})$

This follows by Lemma 1.2.9 from the equality of

$$(-1)^{k-1}e\left(a_{i_1},\dots,a_{i_{k-1}}; A_j, A_l, l \notin \{i_1,\dots,i_{k-1},i,j\}\right) + (-1)^k e\left(a_{i_1},\dots,a_{i_{k-1}},a_i; A_l, l \notin \{i_1,\dots,i_{k-1},i,j\}\right)$$

and

$$(-1)^{k-1} e\big(a_{i_1}, \dots, a_{i_{k-1}}; A_i, A_l, l \notin \{i_1, \dots, i_{k-1}, i, j\}\big) \\ + (-1)^k e\big(a_{i_1}, \dots, a_{i_{k-1}}, a_j; A_l, l \notin \{i_1, \dots, i_{k-1}, i, j\}\big)$$

In fact, their difference is

$$(-1)^{k-1} e\big(a_{i_1}, \dots, a_{i_{k-1}}; a_j A_j, A_l, l \notin \{i_1, \dots, i_{k-1}, i, j\}\big) \\ - (-1)^{k-1} e\big(a_{i_1}, \dots, a_{i_{k-1}}; a_i A_i, A_l, l \notin \{i_1, \dots, i_{k-1}, i, j\}\big)$$

which is zero by the relation $a_1 A_1 + \cdots + a_n A_n = 0$ (and the remark to Lemma 1.2.2).

The general case is obtained by combining these special cases.

For $g \in GL_n(R)$ we put

$$(a^g)^T = g a^T, A^g = A g^{-1},$$

in such a way that $\sum_{i=1}^n a_i^g A_i^g = 0$.

**(3.1.2) Lemma:** If the alternating multiplicities are also defined for $a^g, A^g$ ($g \in GL_n(R)$), then $S(a^g, A^g) = S(a, A)$.

Proof:

(1) If $g$ is an (invertible) diagonal matrix, the assertion is clear (by isomorphism of the defining complexes).

(2) Let $g$ be an elementary matrix, i.e.

$$a_{i_2}^g = a_{i_2} + c a_{i_1}, a_k^g = a_k, k \neq i_2$$

$$A_{i_1}^g = A_{i_1} - c A_{i_2}, A_k^g = A_k, k \neq i_1$$

for some $i_1, i_2 \in \{1, \dots, n\}$, $c \in R$. By Lemma 3.1.1 and the assumption, $S(a^g, A^g, \boldsymbol{i}) = S(a^g, A^g)$ is defined for any sequence $\boldsymbol{i}$. Let $\boldsymbol{i}$ begin with $i_1, i_2$. To show $S(a^g, A^g, \boldsymbol{i}) = S(a, A, \boldsymbol{i})$ we have to compare the summands

$$e\left(a_{i_1}, \dots, a_{i_j}; A_k, k \notin \{i_1, \dots, i_{j+1}\}\right) \text{ and } e\left(a_{i_1}^g, \dots, a_{i_j}^g; A_k^g, k \notin \{i_1, \dots, i_{j+1}\}\right).$$

If $j \geq 2$, both of $a_{i_1}^g, a_{i_2}^g$ appear among $a_{i_1}^g, \dots, a_{i_j}^g$, and the Koszul complex of these elements is isomorphic to the one of $a_{i_1}, \dots, a_{i_j}$. The elements $A_k^g, k \notin \{i_1, \dots, i_{j+1}\}$ are all different from $A_{i_1}^g$. The summands are thus equal. If $j \leq 1$, they are evidently also equal.

(3) For the general case, we may pass to the ring $R^* := R[t]_{\mathfrak{m}[t]}$ which is faithfully flat over $R$.

Assertion 1: Let $\mathfrak{a} \subseteq \mathfrak{m}\tilde{R}$ be an ideal of $\tilde{R} = R[t]_{(\mathfrak{m}, t-r)}$ for some $r \in R$. If the specialization $\mathfrak{a}(r) = \mathfrak{a}\tilde{R}/(t-r) \subseteq R$ is of finite colength, then $\mathfrak{a}R^*$ is of finite colength as well. A similar statement holds for modules $M \subseteq \mathfrak{m}\, \tilde{R}^m$.

Proof: In the opposite case, there is a prime ideal $\mathfrak{p} \subseteq \tilde{R}$ with $\mathfrak{a}R^* \subseteq \mathfrak{p}R^* \subsetneq \mathfrak{m}R^*$. Then $\mathfrak{p} \subsetneq \mathfrak{m}\tilde{R} \subsetneq \mathfrak{m}\tilde{R} + (t-r)\tilde{R}$ is a chain of prime ideals, and the contradiction $\dim \tilde{R}/\mathfrak{p} + (t-r) \geq 1$ follows. The module case follows by consideration of the ideal of $m$-minors.

Now we build the matrix $h = I + t(g - I) \in M_n(R[t])$. Over the local ring $R' \coloneqq R[t]_{(\mathfrak{m},t)}$, $h$ is invertible and admits a decomposition into $h = h_l \dots h_1 h_0$, where $h_0$ is diagonal, $h_1, \dots, h_l$ are elementary, and $h_i mod(t) = I$ $(0 \leq i \leq l)$. By assertion 1 with $r = 0$, over $R^*$ the alternating multiplicities of $a^{h_i \dots h_1 h_0}, A^{h_i \dots h_1 h_0}$, $i = 0, \dots, l$, are defined, since this is the case for $t = 0$. By (1), (2) above

$$S_{R^*}(a^{h_i \dots h_1 h_0}, A^{h_i \dots h_1 h_0}),\ i = 0, \dots, l,$$

is equal to $S_{R^*}(a, A)$. In particular $S_{R^*}(a^h, A^h) = S_{R^*}(a, A)$. Analogously, there is a decomposition $hg^{-1} = \tilde{h}_{\tilde{l}} \dots \tilde{h}_1 \tilde{h}_0$, $\tilde{h}_i mod(t-1) = I$ $(0 \leq i \leq \tilde{l})$, in $R'' \coloneqq R[t]_{(\mathfrak{m},t-1)}$, i.e. $h = \tilde{h}_{\tilde{l}} \dots \tilde{h}_1 \tilde{h}_0 g$. Again, by assertion 1 with $r = 1$, over $R^*$ the alternating multiplicities of $a^{\tilde{h}_i \dots \tilde{h}_1 \tilde{h}_0 g}, A^{\tilde{h}_i \dots \tilde{h}_1 \tilde{h}_0 g}$, $i = 0, \dots, \tilde{l}$, are defined. By (1), (2)

$$S_{R^*}\left(a^{\tilde{h}_i \dots \tilde{h}_1 \tilde{h}_0 g}, A^{\tilde{h}_i \dots \tilde{h}_1 \tilde{h}_0 g}\right) = S_{R^*}(a^g, A^g),\ i = 0, \dots, \tilde{l},$$

in particular $S_{R^*}(a^h, A^h) = S_{R^*}(a^g, A^g)$. Putting this together,

$$S_R(a, A) = S_{R^*}(a, A) = S_{R^*}(a^h, A^h) = S_{R^*}(a^g, A^g) = S_R(a^g, A^g).$$

In the following lemma and its consequences, we assume that $R/\mathfrak{m}$ is infinite by passing to $R[t]_{\mathfrak{m}[t]}$ otherwise.

**(3.1.3) Lemma:** Let $a_1, \dots, a_n \in \mathfrak{m}$ and $A_1, \dots, A_n \in \mathfrak{m}R^m$ with $n = d + m$ elements generating an ideal $\mathfrak{a}$ and a submodule $M$ both of finite colength. For each $1 \leq i \leq d + 1$ there is a matrix $g \in GL_n(R)$ such that

$$e\left(a_1^g, \dots, a_{i-1}^g, A_{i+1}^g, \dots, A_n^g\right) = e(\mathfrak{a}, i - 1; M).$$

Proof: By induction on $i$.

$i = 1$: There is a reduction of $M$ generated by $m + d - 1$ general linear combinations of the generators (cf. appendix). Therefore $g$ may be chosen such that $(A_{i+1}^g, \dots, A_n^g)$ is a reduction of $M$.

$i \to i + 1$: We may assume that $R$ is reduced as follows. Let $\bar{R} = R/\mathfrak{n}$ with $\mathfrak{n}$ the nilradical. Suppose $e\left(a_1^g, \dots, a_{i-1}^g, A_{i+1}^g, \dots, A_n^g; \bar{R}\right) = e(\mathfrak{a}, i - 1; M; \bar{R})$. Then for any prime ideal $\mathfrak{p}$ of dimension $d$, $e\left(a_1^g, \dots, a_{i-1}^g, A_{i+1}^g, \dots, A_n^g; R/\mathfrak{p}\right) = e(\mathfrak{a}, i - 1; M; R/\mathfrak{p})$. Namely, $\geq$ holds by Lemma 1.2.7, and the sums of both sides over all $\mathfrak{p}$ must be equal by Lemma 1.1.3 and Lemma 1.2.6. Then also $e\left(a_1^g, \dots, a_{i-1}^g, A_{i+1}^g, \dots, A_n^g\right) = e(\mathfrak{a}, i - 1; M)$, because the sums with an additional factor $L(R_\mathfrak{p})$ (replacing $L\left(\bar{R}_{\bar{\mathfrak{p}}}\right) = 1$) are equal again.

For the reduced ring $R$ of dimension $d > 0$, the $\mathfrak{m}$ primary ideal $\mathfrak{a}$ does not consist entirely of zero-divisors. Hence we find a $g \in GL_n(R)$ with the properties

(1) $a_1^g$ is a non-zero-divisor,
(2) $a_1^g$ is superficial for $\mathfrak{a}$,$M$ (cf. appendix),
(3) $M' = (A_2^g, \dots, A_n^g)$ is a reduction of $M$ (cf. appendix).

By the induction hypothesis, applied to the ring $\bar{R} = R/a_1^g R$ and the elements $a_2^g, \dots, a_n^g; A_2^g, \dots, A_n^g$, there is a matrix $h \in GL_n(R)$ of the form $h = diag(I_1, \tilde{h})$, $\tilde{h} \in GL_{n-1}(R)$, such that

$$e\big(a_2^{hg},\dots,a_{i-1}^{hg};A_{i+1}^{hg},\dots,A_n^{hg};\bar{R}\big)=e(\mathfrak{a},i-2;M';\bar{R}).$$

Because of (1) and $a_1^g = a_1^{hg}$, the left-hand side is equal to

$$e\big(a_1^{hg},a_2^{hg},\dots,a_{i-1}^{hg};A_{i+1}^{hg},\dots,A_n^{hg};R\big),$$

and by (1),(2),(3) and Lemma 4.2.3, the right hand side is

$$e(\mathfrak{a},i-2;M;\bar{R})=e(\mathfrak{a},i-1;M;R).$$

**(3.1.4) Lemma:** With the same assumptions as in (3.1.3) there is a $g\in GL_n(R^*)\cap M_n(R[t])$, where $R^*=R[t]_{\mathfrak{m}[t]}$, $t=(t_1,\dots,t_{d+1})$, such that

$$e\big(a_1^g,\dots,a_{i-1}^g;A_{i+1}^g,\dots,A_n^g;R^*\big)=e(\mathfrak{a}R^*,i-1;R^*M;R^*)=e(\mathfrak{a},i-1;M;R)$$

for all $1\le i\le d+1$.

Proof: By Lemma 3.1.3 there is such a matrix $g_i$ for each $i$ separately. We take

$$g=\textstyle\sum_{i=1}^{d+1}t_ig_i\in GL_n(R^*)\cap M_n(R[t]).$$

Assertion 1: Let $b_1,\dots,b_{k-1}\in\mathfrak{m}[t]$, $B_{k+1},\dots,B_{m+d}\in\mathfrak{m}R[t]^m$, such that $e(b,B;R^*)$ is defined. Then

$$e(b,B;R^*)\le e(t-r,b,B;R')=e(b(r),B(r);R)$$

for all $r\in R^{d+1}$ such that the last term is defined. Here $R'=R[t]_{(\mathfrak{m},t-r)}$, and $b(r),B(r)$ are defined by the images in $R[t]/(t-r)\cong R$.

Proof: Let $L$ be the complex defined by $b,B$ over $R'$ and $K$ the Koszul complex of $s:=t-r$. Then with $\mathfrak{p}=\mathfrak{m}[t]R'$,

$$e(s,b,B;R')=\chi(K\otimes L)=\chi(K\otimes H(L))\ge\textstyle\sum_i(-1)^i\chi(K\otimes R'/\mathfrak{p})L(H^i(L)_\mathfrak{p})=\chi\big(L\otimes R'_\mathfrak{p}\big)=\chi(L\otimes R^*)=e(b,B;R^*),$$

where we have used Lemma 1.2.6 and $R'/\mathfrak{p}\cong R/\mathfrak{m}[t]_{(t-r(0))}$ (cf. Lemma 1.2.11). The inequality $\ge$ holds, because we only take care of one of the $(d+1)$-dimensional prime ideals in the support of $H(L)$ (of dimension $d+1$) and the left-out contributions are nonnegative (e.g. by Lemma 1.2.7). The equality $e(s,b,B;R')=e(b(r),B(r);R)$ follows from Lemma 1.2.8, since $s$ is a regular sequence in $R'$.

From assertion 1 we obtain (by specializing $t_j=\delta_{ij}$)

$$e\big(a_1^g,\dots,a_{i-1}^g;A_{i+1}^g,\dots,A_n^g;R^*\big)\le e\big(a_1^{g_i},\dots,a_{i-1}^{g_i};A_{i+1}^{g_i},\dots,A_n^{g_i};R\big)=e(\mathfrak{a},i-1,M)=e(\mathfrak{a}R^*,i-1,R^*M;R^*).$$

By Lemma 1.2.7 there is equality.

### 3.2 Final step of the proof theorem 0.2

Let $\mathfrak{a}=(a_1,\dots,a_n)$, $M=(A_1,\dots,A_n)$. By Lemma 3.1.4 there is a $g\in GL_n(R^*)$ with the equality

$$e(\mathfrak{a}, i-1; M; R) = e\left(a_1^g, \dots, a_{i-1}^g; A_{i+1}^g, \dots, A_n^g; R^*\right), 1 \le i \le d+1.$$

By Theorem 0.1, together with Lemma 3.1.2 and the faithful flatness of $R^*$ over $R$ we conclude

$$\chi(K_\nu) = S_{R^*}(a^g, A^g) = S_{R^*}(a, A) = S(a, A).$$

## 4 Appendix: Reductions · Superficial elements

Although these matters are well-known (cf. e.g. [Ga], [KT] and the literature cited) we give a complete justification of the facts used in the main part.

### 4.1 Reductions

Let $R, \mathfrak{m}$ be a Noetherian local ring and $M = (A_1, \dots, A_n) \subseteq \mathfrak{m}R^m$ a module of finite colength. We consider the generated algebra

$$U = \oplus_{\nu \ge 0} U_\nu = R[y_1, \dots, y_n] \subseteq S = R[x_1, \dots, x_m],$$

$$y_j = \sum_{i=1}^m A_{ij}\, x_i, j = 1, \dots, n.$$

**(4.1.1) Lemma:** The minimal prime ideals of $U$ are all $\mathfrak{p}^* = \mathfrak{p}S \cap U$, $\mathfrak{p}$ a minimal prime ideal of $R$. The quotient ring $U/\mathfrak{p}^*$ is generated by $\bar{M} = \bar{R}M \subseteq \bar{R}^m$ in $\bar{R}[x_1, \dots, x_m]$, where $\bar{R} = R/\mathfrak{p}$.

Proof: Let $\mathfrak{P} \subseteq U$ be minimal and $\mathfrak{p} = \mathfrak{P} \cap U_0 \subseteq R$. Then $\mathfrak{P} = Ann(u)$ for some homogeneous element $u \in U$. We have the implications $\mathfrak{p}u = 0 \Rightarrow \mathfrak{p}^*u = 0 \Rightarrow \mathfrak{p}^* \subseteq \mathfrak{P} \Rightarrow \mathfrak{p}^* = \mathfrak{P}$. The description of $U/\mathfrak{p}^*$ is now clear.

**(4.1.2) Lemma:** $\dim U \le d + m$ and $\dim U \ge d + m$ if $d > 0$.

Proof: By Lemma 4.1.1 we may assume $R$ integral of dimension $d > 0$. Over the quotient field $Q = Q(R)$ the elements $y_1, \dots, y_n$ generate $Q^m$, and $Q[y_1, \dots, y_n] = Q[x_1, \dots, x_m]$. Thus $\mathrm{trdeg}_R U = m$. By the dimension formula [ZS, app. 1] for any prime ideal $\mathfrak{P} \subseteq U$ and $\mathfrak{p} = \mathfrak{P} \cap U_0$,

$$ht\, \mathfrak{P} + \mathrm{trdeg}_{R/\mathfrak{p}} U/\mathfrak{P} \le ht\, \mathfrak{p} + \mathrm{trdeg}_R U,$$

in particular $ht\, \mathfrak{P} \le d + m$.

For the converse inequality we let $\mathfrak{P}_i = U \cap (x_1, \dots, x_i)S \subseteq U$. The quotient ring of this prime ideal is

$$U/\mathfrak{P}_i \cong U + (x_1, \dots, x_i)/(x_1, \dots, x_i)S = R[\bar{y}_1, \dots, \bar{y}_n] \subseteq R[x_{i+1}, \dots, x_m],$$

$$\bar{y}_j = \sum_{k=i+1}^m A_{kj} x_k,$$

i.e. we obtain the ring generated by $(M + Re_1 + \dots + Re_i)/Re_1 + \dots + Re_i \subseteq R^{m-i}$. By induction on $m$ (with trivial beginning $m = 0$) we can assume that $\dim U/\mathfrak{P}_1 \ge d + m - 1$. Since $R^m/M$ is of finite colength, there is a nonzero $a \in R$ with $aR^m \subseteq M$. Then $0 \ne ax_1 \in \mathfrak{P}_1$. By our assumption that $R$ is an integral domain, $\dim U \ge d + m$.

**(4.1.3) Lemma:** $\dim U/\mathfrak{m}U \le d + m - 1$, for $d > 0$ there is equality.

Proof: $\leq$ :In the case $d > 0$, for $\mathfrak{p}$ a minimal prime, $\mathfrak{m}U \not\subseteq \mathfrak{p}^*$, hence $\dim U/\mathfrak{m}U < \dim U$. If $d = 0$, $\dim U = \dim U/\mathfrak{m}U = 0$ as follows from $U_\nu = 0, \nu \gg 0$, and $L(U) < \infty$.

$\geq$ : It is sufficient to treat the case $R/\mathfrak{m}$ infinite as follows: Let $\tilde{R} = R[t]_{\mathfrak{m}[t]}$ with maximal ideal $\tilde{\mathfrak{m}} = \mathfrak{m}\tilde{R}$, and $\tilde{U} = U \otimes_R \tilde{R} \subseteq \tilde{R}[x]$. Then

$$\tilde{U}/\tilde{\mathfrak{m}}\tilde{U} = U/\mathfrak{m}U \otimes \tilde{R} = (U/\mathfrak{m}U) \otimes_{R/\mathfrak{m}} (R/\mathfrak{m} \otimes \tilde{R}) = (U/\mathfrak{m}U) \otimes_{R/\mathfrak{m}} (\tilde{R}/\tilde{\mathfrak{m}}),$$

and $\dim \tilde{U}/\tilde{\mathfrak{m}}\tilde{U} = \dim U/\mathfrak{m}U$ from the Hilbert function.

By the normalization theorem [ZS, vol. I, ch. V, th. 8] there are $l := \dim U/\mathfrak{m}U$ elements $u_1, \dots, u_l \in U_1$, whose residue classes generate a subring $R/\mathfrak{m}[\bar{u}_1, \dots, \bar{u}_l] \subseteq U/\mathfrak{m}U$ such that $U/\mathfrak{m}U$ is finite over this subring, i. e. $U_\nu \subseteq (u_1, \dots, u_l)U, \nu \gg 0$. Then the submodule $N := (u_1, \dots, u_l) \subseteq M \cong U_1$ in $R^m$ has finite colength. In fact, otherwise there is a prime ideal $\mathfrak{p} \subsetneq \mathfrak{m}$ such that $S_1/N \otimes k(\mathfrak{p}) \neq 0$. Then also $S_\nu/NS_{\nu-1} \otimes k(\mathfrak{p}) \neq 0$, contradicting $S_\nu/U_\nu \otimes k(\mathfrak{p}) = 0$ and $U_\nu = NU_{\nu-1}, \nu \gg 0$. However, a submodule of finite colength has $n_1 \geq d + m - 1$ generators (if $d > 0$) since $n_1 \geq m$ by rank and the ideal of $m$-minors has height $\leq n_1 - m + 1$.

A *reduction* of $M$ is by definition a submodule $N \subseteq U_1 \cong M$, which generates in $U$ an irrelevant ideal, i.e. $NU_\nu = U_{\nu+1}, \nu \gg 0$. We just have shown:

**(4.1.4) Lemma:** If $R/\mathfrak{m}$ is infinite and $d > 0$, then there is a reduction $N$ of $M$, which is generated by $d + m - 1$ elements. Such a reduction is minimal.

We note that by the proof of the normalization theorem, $m + d - 1$ “general” linear combinations of generators of $M$ generate a reduction. By “general” one means that the coefficient matrix reduced modulo $\mathfrak{m}$ is allowed to lie in a suitable non-empty open set.

### 4.2 Superficial elements

For submodules $\mathfrak{a} \subseteq \mathfrak{m}$, $M \subseteq \mathfrak{m}R^m$ of finite colength over the local ring $R$ of dimension $d$ we keep the definition of the rings $S$ and $U$ as in 4.1 and newly form the rings

$$\hat{S} := \oplus_{\mu,\nu \geq 0} s^\mu S_\nu = R[x, s],\ \hat{U} := \oplus_{\mu,\nu \geq 0} \mathfrak{a}^\mu s^\mu U_\nu \subseteq R[x, s],$$

$$\bar{U} := \oplus_{\mu,\nu \geq 0} \mathfrak{a}^\mu U_\nu / \mathfrak{a}^{\mu+1} U_\nu.$$

The element $a \in \mathfrak{a}$ is called *superficial* for $\mathfrak{a}, M$, if the multiplication map by $\bar{a} \in \mathfrak{a}/\mathfrak{a}^2$,

$$\bar{U}_{\mu\nu} \to \bar{U}_{\mu+1,\nu},\ \mu, \nu \gg 0,$$

is injective. This means:

$$\exists \mu_0, \nu_0 \colon (\mathfrak{a}^{\mu+2} U_\nu \colon a) \cap \mathfrak{a}^{\mu_0} U_\nu = \mathfrak{a}^{\mu+1} U_\nu,\ \mu \geq \mu_0, \nu \geq \nu_0.$$

The condition is satisfied in particular, if $\bar{a}$ is a superficial element for the ring $\bar{U}$, which we view as single-graded by $\mu$, i.e. if $\bar{a}$ is for sufficiently large $\mu$ (and arbitrary $\nu$) a non-zero-divisor.

This sharper condition is equivalent to the requirement that $\bar{a}$ should not lie in the union of all relevant associated prime ideals of $\bar{U}$ (i.e. those not containing $\bar{U}_+$). If $R/\mathfrak{m}$ is infinite, the union of finitely many proper subspaces of $\mathfrak{a}/\mathfrak{m}\mathfrak{a}$ is a proper subset, and the requirement is fulfilled for “general” elements.

The following two lemmas are mixed analogues of the usual theorems for ordinary multiplicities.

**(4.2.1) Lemma:** Let $a \in \mathfrak{a}$ be a superficial element for $\mathfrak{a}, M$, and assume $d \geq 2$. Then $\bar{R} = R/aR$ is of dimension $d-1$, and

$$e(\mathfrak{a}, i, M) = e(\bar{\mathfrak{a}}, i-1, \bar{M}),\ 2 \leq i \leq d.$$

Proof: By the assumption, for $\mu, \nu \gg 0$ the sequence

$$0 \to \mathfrak{a}^{\mu}U_{\nu}/\mathfrak{a}^{\mu+1}U_{\nu} \overset{\bar{a}}{\to} \mathfrak{a}^{\mu+1}U_{\nu}/\mathfrak{a}^{\mu+2}U_{\nu} \to \mathfrak{a}^{\mu+1}U_{\nu}/\mathfrak{a}^{\mu+2}U_{\nu} + a\mathfrak{a}^{\mu}U_{\nu} \to 0$$

is exact. The task is to rewrite the right term suitably.

Assertion 1: $(\mathfrak{a}^{\mu+1}U_{\nu}: a)_{S_{\nu}} = (0: a)_{S_{\nu}} + \mathfrak{a}^{\mu}U_{\nu},\ \mu, \nu \gg 0.$

Proof: By the fact, that $a\hat{S} \cap \widehat{U}$ is a finitely generated ideal of $\widehat{U}$, there are $\mu_0, \nu_0$ such that

$$aS_{\nu} \cap \mathfrak{a}^{\mu+1}U_{\nu} = \mathfrak{a}^{\mu-\mu_0+1}U_{\nu-\nu_0}(aS_{\nu_0} \cap \mathfrak{a}^{\mu_0}U_{\nu_0}),\ \mu \geq \mu_0, \nu \geq \nu_0.$$

Let $z \in S_{\nu}$ with $az \in \mathfrak{a}^{\mu+1}U_{\nu}$ an element of $(\mathfrak{a}^{\mu+1}U_{\nu}: a)_{S_{\nu}}$. Then $az = au,\ u \in \mathfrak{a}^{\mu-\mu_0+1}U_{\nu-\nu_0}S_{\nu_0}$, hence $a(z-u) = 0,\ z \equiv u \ mod\ (0: a)_{S_{\nu}}$. If $\mu_0, \nu_0$ are large enough,

$$(\mathfrak{a}^{\mu+1}U_{\nu}: a)_{U_{\nu}} \cap \mathfrak{a}^{\mu_0}U_{\nu} = \mathfrak{a}^{\mu}U_{\nu},\ \mu \geq \mu_0, \nu \geq \nu_0.$$

Since $S_{\nu_0}/U_{\nu_0}$ is of finite length, $u \in \mathfrak{a}^{\mu-\mu_0+1}U_{\nu-\nu_0}S_{\nu_0} \subseteq \mathfrak{a}^{\mu_0}U_{\nu}$ for $\mu \gg 0$, and then as required $u \in (\mathfrak{a}^{\mu+1}U_{\nu}: a)_{U_{\nu}} \cap \mathfrak{a}^{\mu_0}U_{\nu} = \mathfrak{a}^{\mu}U_{\nu}$.

By assertion 1, for $\mu, \nu \gg 0$ there is equality

$$\mathfrak{a}^{\mu+1}U_{\nu}/\mathfrak{a}^{\mu+2}U_{\nu} + a\mathfrak{a}^{\mu}U_{\nu} \overset{=}{\to} \mathfrak{a}^{\mu+1}U_{\nu}/\mathfrak{a}^{\mu+2}U_{\nu} + aS_{\nu} \cap \mathfrak{a}^{\mu+1}U_{\nu},$$

which gives the exact sequence

$$0 \to \bar{U}_{\mu\nu} \to \bar{U}_{\mu+1,\nu} \to \bar{U}_{\bar{R},\mu+1,\nu} \to 0,\ \mu, \nu \gg 0,$$

where $\bar{U}_{\bar{R}}$ denotes the corresponding ring for $\bar{R}$ and the images $\bar{\mathfrak{a}}, \bar{M}$. For the Hilbert polynomials one deduces

$$\Delta_{\mu}^{2}Q(\mu-1, \nu) = \Delta_{\mu}\bar{Q}(\mu, \nu),$$

and from the leading terms one reads off the claim.

**(4.2.2) Lemma:** In the situation of Lemma 4.2.1, if $a$ additionally is a non-zero-divisor, then for all $d \geq 1$ and $1 \leq i \leq d$

$$e(\mathfrak{a}, i, M) = e(\bar{\mathfrak{a}}, i-1, \bar{M}).$$

Proof: We contemplate the exact sequence

$$0 \to (\mathfrak{a}^{\mu+1}U_{\nu}: a)_{S_{\nu}}/\mathfrak{a}^{\mu}U_{\nu} \to S_{\nu}/\mathfrak{a}^{\mu}U_{\nu} \overset{\bar{a}}{\to} S_{\nu}/\mathfrak{a}^{\mu+1}U_{\nu} \to S_{\nu}/aS_{\nu} + \mathfrak{a}^{\mu+1}U_{\nu} \to 0.$$

By assertion 1 of the last proof the left term is zero for $\mu, \nu \gg 0$. We obtain $\Delta_{\mu}Q(\mu, \nu) = \bar{Q}(\mu+1, \nu)$.

We remark, that the exact sequence provides an inequality of the multiplicities if $a$ is merely a dimension lowering element.

**(4.2.3) Lemma:** Let $M'$ be a reduction of $M$. Then $e(\mathfrak{a}, i, M') = e(\mathfrak{a}, i, M)$.

Proof: With the help of the last two lemmas, this can be reduced to the unmixed case, which we prove apart.

**(4.2.4) Lemma:** Let $M'$ be a reduction of $M$. Then $e(M') = e(M)$.

Proof: Writing $N = M'$, we have to show that the difference $L(M^\nu/N^\nu) = L(S_\nu/N^\nu) - L(S_\nu/M^\nu)$, $\nu \gg 0$, is a polynomial of degree $< d + m - 1$. There is a constant $\nu_0$ such that $M^{\nu_0+\nu} = M^{\nu_0}N^\nu$, $\nu \geq 0$. Over the ring $U' := \bigoplus_{\nu\geq 0} N^\nu$ we consider the finitely generated module

$$\bigoplus_{\nu\geq 0} {}^{M^{\nu_0}N^\nu}\!/_{N^{\nu_0}N^\nu}.$$

Since $M^{\nu_0}/N^{\nu_0}$ is of finite length, we have a module over $U'/\mathfrak{m}^k U'$ for some $k > 0$. By Lemma 4.1.3, $\dim U'/\mathfrak{m}U' \leq d + m - 1$, and the claim follows.